\documentclass[1p]{elsarticle}
\usepackage{amssymb}
\usepackage{amsfonts}

\newtheorem{theorem}{Theorem}
\newtheorem{conjecture}[theorem]{Conjecture}

\newtheorem{lemma}[theorem]{Lemma}

\newtheorem{remark}{Remark}

\newproof{pf}{Proof}

\begin{document}
\title{Distant irregularity strength of graphs with bounded minimum degree}

\author{Jakub Przyby{\l}o\fnref{fn1,fn2}}
\ead{jakubprz@agh.edu.pl, phone: 048-12-617-46-38,  fax: 048-12-617-31-65}

\fntext[fn1]{Financed within the program of the Polish Minister of Science and Higher Education
named ``Iuventus Plus'' in years 2015-2017, project no. IP2014 038873.}
\fntext[fn2]{Partly supported by the Polish Ministry of Science and Higher Education.}

\address{AGH University of Science and Technology, al. A. Mickiewicza 30, 30-059 Krakow, Poland}

\begin{abstract}
Consider a graph $G=(V,E)$ without isolated edges and with maximum degree $\Delta$.
Given a colouring $c:E\to\{1,2,\ldots,k\}$, the weighted degree of a vertex $v\in V$
is the sum of its incident colours, i.e., $\sum_{e\ni v}c(e)$.
For any integer $r\geq 2$,
the least $k$ admitting the existence of such $c$ attributing distinct weighted degrees to
any two different vertices at distance at most $r$ in $G$ is called the $r$-distant irregularity strength of $G$ and denoted by $s_r(G)$.
This graph invariant provides a natural link between the well known 1--2--3 Conjecture and irregularity strength of graphs.
In this paper we apply the probabilistic method in order to prove an upper bound $s_r(G)\leq (4+o(1))\Delta^{r-1}$ for graphs with minimum degree $\delta\geq \ln^8\Delta$, improving thus far best upper bound $s_r(G)\leq 6\Delta^{r-1}$.
\end{abstract}


\maketitle

\section{Introduction}
Let us consider a graph $G=(V,E)$ and its not necessarily proper edge colouring $c:E\to\{1,2,\ldots,k\}$ with $k$ least positive integers.
We say that such colouring $c$ is \emph{irregular} if it associates with every vertex $v\in V$ a different
sum of its incident colours:
\begin{equation}\label{DefOfW}
w_c(v):=\sum_{u\in N(v)}c(uv),
\end{equation}
so called \emph{weighted degree} of $v$. We shall also denote $w_c(v)$ by $w(v)$ 
in cases when the colouring $c$ is unambiguous from context.
The least $k$ admitting such irregular colouring $c$ is called the \emph{irregularity strength} of $G$ and denoted by $s(G)$, see~\cite{Chartrand}.
Note that this parameter is well defined for graphs without isolated edges and with at most one isolated vertex; for the remaining ones we might e.g. set $s(G)=\infty$.
Alternatively, $s(G)$ might be regarded as the least $k$ so that we may construct an irregular multigraph, i.e. a multigraph with pairwise distinct degrees of all vertices, of $G$ by multiplying its edges, each at most $k$ times.
This study thus originate from the basic fact that no graph $G$ with more than one vertex is irregular itself, hence $s(G)\geq 2$, and related research on possible alternative definitions of irregularity in graph environment, see e.g.~\cite{ChartrandErdosOellermann}.
It is known that $s(G)\leq n-1$, where $n=|V|$, for all graphs containing no isolated edges and at most one isolated vertex, except for the graph $K_3$, see \cite{Aigner,Nierhoff}. This is a tight upper bound, as exemplified e.g. by the family of stars. A better upper bound is known for graphs with minimum degree $\delta>6$, i.e., $s(G)\leq 6\lceil\frac{n}{\delta}\rceil$, see
\cite{KalKarPf}, and $s(G)\leq (4+o(1))\frac{n}{\delta}+4$ for graphs with
$\delta\geq n^{0.5}\ln n$, see \cite{MajerskiPrzybylo2}.
It is however believed that these upper bounds can be improved to (in such a case optimal) $s(G)\leq \frac{n}{\delta} + C$ for some absolute constant $C$, see e.g. \cite{Faudree,KalKarPf,MajerskiPrzybylo2,irreg_str2}. This has been explicitly conjectured in the case of $d$-regular graphs, see \cite{Faudree}, for which one can
 observe that on the other hand $s(G)\geq \frac{n}{d}+\frac{d-1}{d}$ via straightforward counting argument, see e.g.~\cite{Chartrand}.
Other results concerning the concept of irregularity strength and in particular its value for specific graph classes can also be found e.g. in \cite{Bohman_Kravitz,Lazebnik,Dinitz,Faudree2,Frieze,Lehel,Przybylo}, and many others.

The problem described above gave rise to a variety of associated questions and concepts, nowadays making up an intensively studied field of the graph theory. One of its most intriguing descendants is its local version, where one investigates the least $k$ for which
there is a colouring $c:E\to\{1,2,\ldots,k\}$ of the edges of a graph $G$ such that $w_c(u)\neq w_c(v)$ for every pair of adjacent vertices $u,v$, so called \emph{neighbours} in $G$. Denote this value by $s_1(G)$ (we shall comment on this notion below), and note that such graph invariant is well defined for all graphs without isolated edges. Though initially no finite upper bound was known for this parameter, Karo\'nski, {\L}uczak and Thomason~\cite{123KLT} posed a fascinating conjecture that $s_1(G)\leq 3$ for all graphs without isolated edges. This is nowadays commonly referred to as \emph{1--2--3 Conjecture} in the literature, see e.g.~\cite{KalKarPf_123}.
The conjecture is still open,
while thus far the following general upper bounds were subsequently proved: $s_1(G)\leq 30$ in~\cite{Louigi30},
$s_1(G)\leq 16$ in~\cite{Louigi}, $s_1(G)\leq 13$ in~\cite{123with13},
and finally $s_1(G)\leq 5$ from~\cite{KalKarPf_123}.

In this paper we study a problem linking the two concepts above.
Given any graph $G=(V,E)$ and an integer $r\geq 1$, two distinct vertices at distance at most $r$ in $G$, i.e. $u,v\in V$ with $1\leq d(u,v)\leq r$, shall be called \emph{$r$-neighbours}. For any colouring $c:E\to\{1,2,\ldots,k\}$, the weighted degree $w_c(v)$ (defined in~(\ref{DefOfW})) shall also be referred to as
the \emph{weight} of $v$ or simply the \emph{sum} at $v$. If $w(u)=w(v)$ for distinct vertices $u,v\in V$, we say that they are \emph{in conflict}, otherwise we call them \emph{sum-distinguished} or simply \emph{distinguished}. The least $k$ such that there is a colouring $c:E\to\{1,2,\ldots,k\}$ without a conflict between any pair of $r$-neighbours in the graph $G$ shall be called the \emph{$r$-distant irregularity strength} of $G$ and denoted by $s_r(G)$ (observe that this notion is consistent with the use of $s_1(G)$ 
 with reference to 1--2--3 Conjecture above, while it is also justified to set $s_\infty(G)=s(G)$ in this context). Note that analogously as above, $s_r(G)$ is well defined iff $G$ has no isolated edges.
In~\cite{Przybylo_distant} the following upper bound was provided for this graph invariant.
\begin{theorem}\label{Old6Theorem}
Let $G$ be a graph without isolated edges, and with maximum degree $\Delta\geq 2$, and let $r\geq 1$ be an integer. Then,
$$s_r(G)\leq 6\Delta^{r-1}.$$
\end{theorem}
See also~\cite{Przybylo_distant} for a discussion justifying the fact that the general upper bound from the theorem above cannot be smaller than $\Delta^{r-1}$.
In this paper we essentially improve the inequality from Theorem~\ref{Old6Theorem} to $s_r(G)\leq (4+o(1))\Delta^{r-1}$, but for technical reasons
we had to exclude from our result graphs with
small minimum degrees, i.e., smaller than a value given by a certain poly-logarithmic function in $\Delta$,
see Theorem~\ref{przybylo_main_th} below.
Our approach shall be based on the probabilistic method, first applied to design a special ordering of the vertices of a graph,
and then to provide an enhancement of
an algorithm whose different variants were used e.g. in~\cite{KalKarPf,KalKarPf_123,MajerskiPrzybylo2,Przybylo_distant},
developed along the specified order.
In the next section we recall several useful tools of the probabilistic method.
Then we formulate our main result, and provide its proof in Section~\ref{SectionWithPrzybyloProof}. The last section contains a few related comments.

\section{Tools}
We shall use a few tools of the probabilistic method listed in details below.
In particular,
the Lov\'asz Local Lemma, see e.g.~\cite{AlonSpencer},
combined with the Chernoff Bound, see e.g.~\cite{JansonLuczakRucinski}
(Th. 2.1, page 26)
and Talagrand's Inequality, see e.g.~\cite{MolloyReed_GoodTalagrand}. 
\begin{theorem}[\textbf{The Local Lemma}]
\label{LLL-symmetric}
Let $A_1,A_2,\ldots,A_n$ be events in an arbitrary pro\-ba\-bi\-li\-ty space.
Suppose that each event $A_i$ is mutually independent of a set of all the other
events $A_j$ but at most $D$, and that ${\rm \emph{\textbf{Pr}}}(A_i)\leq p$ for all $1\leq i \leq n$. If
$$ ep(D+1) \leq 1,$$
then $ {\rm \emph{\textbf{Pr}}}\left(\bigcap_{i=1}^n\overline{A_i}\right)>0$.
\end{theorem}
\begin{theorem}[\textbf{Chernoff Bound}]\label{ChernofBoundTh}
For any $0\leq t\leq np$,
$${\rm\emph{\textbf{Pr}}}({\rm BIN}(n,p)>np+t)<e^{-\frac{t^2}{3np}}~~{and}~~{\rm\emph{\textbf{Pr}}}({\rm BIN}(n,p)<np-t)<e^{-\frac{t^2}{2np}}\leq e^{-\frac{t^2}{3np}}$$
where ${\rm BIN}(n,p)$ is the sum of $n$ independent Bernoulli variables, each equal to $1$ with probability $p$ and $0$ otherwise.
\end{theorem}

\begin{theorem}[\textbf{Talagrand's Inequality}]\label{TalagrandsInequalityTotal}
Let $X$ be a non-negative random variable determined by $l$ independent trials $T_1,\ldots,T_l$. 
Suppose there exist constants $c,k>0$ such that for every set of possible outcomes of the trials, we have:
\begin{itemize}
\item[1.] changing the outcome of any one trial can affect $X$ by at most $c$, and
\item[2.] for each $s>0$, if $X\geq s$ then there is a set of at most $ks$ trials whose outcomes certify that $X\geq s$. 
\end{itemize}
Then for any $t\geq 0$, we have
$${\rm\emph{\textbf{Pr}}}(|X-{\rm\emph{\textbf{E}}}(X)|>t+20c\sqrt{k{\rm\emph{\textbf{E}}}(X)}+64c^2k)\leq 4e^{-\frac{t^2}{8c^2k(\mathbf{E}(X)+t)}}.$$
\end{theorem}

Note that e.g. knowing only an upper bound $\mathbf{E}(X)\leq h$ (instead of the exact value of $\mathbf{E}(X)$)
we may still use Talagrand's Inequality in order to upper-bound the probability that $X$ is large.
It is sufficient to apply Theorem~\ref{TalagrandsInequalityTotal} above
to the variable $Y=X+h-\mathbf{E}(X)$, with $\mathbf{E}(Y)=h$ to obtain the following
provided that the assumptions of Theorem~\ref{TalagrandsInequalityTotal}
hold for $X$: 
$$\mathbf{Pr}(X>h+t+20c\sqrt{kh}+64c^2k) \leq \mathbf{Pr}(Y>h+t+20c\sqrt{kh}+64c^2k) \leq 4e^{-\frac{t^2}{8c^2k(h+t)}}.$$
Analogously, in the case of the Chernoff Bound, if $X$ is a sum of $n\leq k$ (where $k$ does not have to be an integer)
random independent Bernoulli variables, each equal to $1$ with probability $p\leq q$,
then $\mathbf{Pr}(X>kq+t)< e^{-\frac{t^2}{3kq}}$ (for $t\leq\lfloor k\rfloor q$).

\section{Main Result}
\begin{theorem}\label{przybylo_main_th}
For every integer $r\geq 2$ there exists a constant $\Delta_0$ such that for each
graph $G$ with maximum degree $\Delta\geq \Delta_0$ and minimum degree $\delta\geq \ln^8\Delta$,
$$s_r(G) < 4\Delta^{r-1}\left(1+ \frac{1}{\ln\Delta}\right) + 12,$$
hence $s_r(G)\leq (4+o(1))\Delta^{r-1}$ for all graphs with $\delta\geq \ln^8\Delta$ and without isolated edges.
\end{theorem}

We do not specify the value of $\Delta_0$ in
the proof below (nor in the statement of the theorem above),
assuming whenever needed that $\Delta$ is large enough so that some explicit inequalities hold.
Note also that the conclusion from the last line 
of Theorem~\ref{przybylo_main_th} above follows by the upper bound from Theorem~\ref{Old6Theorem} (applied to graphs with $\Delta<\Delta_0$).

\section{Proof of Theorem~\ref{przybylo_main_th}\label{SectionWithPrzybyloProof}}
Fix any integer $r\geq 2$ and let $G=(V,E)$ be a graph
with minimum degree $\delta\geq \ln^8\Delta$, where $\Delta$ is the (sufficiently large) maximum degree of $G$.

Let $Q$ and $q$ be the least integers divisible by $3$ such that
$Q \geq 2\Delta^{r-1}+ \frac{\Delta^{r-1}}{\ln\Delta}$ and $q \geq \frac{\Delta^{r-1}}{\ln\Delta}$
(hence $Q < 2\Delta^{r-1}+ \frac{\Delta^{r-1}}{\ln\Delta}+3$ and $q < \frac{\Delta^{r-1}}{\ln\Delta}+3$).
We shall show that $s_r(G) \leq 2Q+2q$.

Let $V=\{v_1,v_2,\ldots,v_n\}$ (where $v_i\neq v_j$ for $i\neq j$).
First we shall randomly (re)order the vertices in $V$.
For this goal, independently for every $i=1,\ldots,n$, we pick a (real) number uniformly at random from the interval $[0,1]$
and associate it with $v_i$.
As a result we obtain a (new) ordering $u_1,\ldots,u_n$ where a vertex $u_i$ precedes $u_j$ if and only if the value chosen for $u_i$ is not greater than the one chosen for $u_j$.
In other words, we associate with every vertex $v$ a random variable $X_v\sim U[0,1]$ having the uniform distribution on $[0,1]$,
and order the vertices in $V$ into a sequence $u_1,\ldots,u_n$ so that
$X_{u_i}\leq X_{u_j}$ whenever $i\leq j$ (i.e.,
subsequent vertices in this ordering correspond to
the order statistics of the defined set of independent random variables).
Note also that we may assume that such ordering is uniquely defined, i.e. that $X_u\neq X_v$ whenever $u\neq v$ with probability one
(as the probability that $X_v=X_{v'}$ for some pair $v,v'$ of distinct vertices in $V$ equals $0$).

Let us partition $V$ into three subsets:
\begin{eqnarray}
A&=&\left\{v:X_v<\frac{1}{\ln^2\Delta}\right\};\nonumber\\
B&=&\left\{v:\frac{1}{\ln^2\Delta}\leq X_v\leq 1-\frac{1}{\ln^{3}\Delta}\right\};\nonumber\\
C&=&\left\{v:X_v> 1-\frac{1}{\ln^{3}\Delta}\right\}.\nonumber
\end{eqnarray}

For every vertex $v\in V$, any its neighbour or $r$-neighbour $u$ which precedes $v$ in the obtained 
ordering of the elements of $V$ shall be called a \emph{backward neighbour} or \emph{$r$-neighbour}, resp., of $v$. Analogously, the remaining ones shall be called \emph{forward neighbours} or \emph{$r$-neighbours}, resp., of $v$, while the edges joining $v$ with its forward or backward neighbours shall be referred to as \emph{forward} or \emph{backward}, resp., as well.
Also, for any subset $S\subset V$, let $N_-(v)$, $N_-^r(v)$, $N_S(v)$, $N_S^r(v)$ denote the sets of all backward neighbours, backward $r$-neighbours, neighbours in $S$ and $r$-neighbours in $S$ of $v$, respectively.
Set
$d_-(v)=|N_-(v)|$, $d_-^r(v)=|N_-^r(v)|$, $d_S(v)=|N_S(v)|$, $d_S^r(v)=|N_S^r(v)|$, and for any subset of edges $E_0\subseteq E$, set $d_{E_0}(v)=|\{u\in N(v):uv\in E_0\}|$.\\
\begin{lemma}\label{MainSequencingLemma}
With positive probability, the obtained ordering has the following features for every vertex $v$ in $G$, whose degree we denote by $d$:
\begin{itemize}
\item[$F_1$:]
$d^r_A(v) \leq 2\frac{d\Delta^{r-1}}{\ln^2\Delta}$;
\item[$F_2$:]
$d^r_C(v)\leq 2\frac{d\Delta^{r-1}}{\ln^{3}\Delta}$;
\item[$F_3$:]
$\frac{1}{2}\frac{d}{\ln^2\Delta}
\leq d_A(v) \leq 2\frac{d}{\ln^2\Delta}$;
\item[$F_4$:]
$\frac{1}{2}\frac{d}{\ln^{3}\Delta}
\leq d_C(v) \leq 2\frac{d}{\ln^{3}\Delta}$;
\item[$F_5$:] if $v\in B$, then: $d_-(v)\geq X_v d-\sqrt{X_v d}\ln\Delta$;
\item[$F_6$:] if $v\in B$, then: $d^r_-(v)\leq X_v d \Delta^{r-1}+\sqrt{X_vd\Delta^{r-1}}\ln\Delta$.
\end{itemize}
\end{lemma}

\begin{pf} For every vertex $v\in V$ of degree $d\in[\delta,\Delta]$, let $A_{v,1},A_{v,2},A_{v,3},A_{v,4}$ denote the events that the features $F_1,F_2,F_3,F_4$, resp., are not satisfied for $v$. Denote by $A_{v,5}$ the event that $v$ belongs to $B$ and $d_-(v) < X_v d-\sqrt{X_v d}\ln\Delta$, and analogously,
let $A_{v,6}$ denote the event that $v$ belongs to $B$ and $d^r_-(v) > X_v d \Delta^{r-1}+\sqrt{X_vd\Delta^{r-1}}\ln\Delta$.
Exploiting the Lov\'asz Local Lemma we
shall exhibit that with positive probability none of the events $A_{v,i}$ holds for any $v\in V$ and $i\in\{1,\ldots,6\}$,
thus proving the thesis of Lemma~\ref{MainSequencingLemma}.

Note that for every vertex $v\in V$ of degree $d$, $d^r(v)\leq d\Delta^{r-1}$. Moreover, for each $u\in N^r(v)$, the probability that $u$ belongs to $A$ equals $\frac{1}{\ln^2\Delta}$. Thus by the Chernoff Bound (and the comments below it),
\begin{equation}\mathbf{Pr}(A_{v,1}) \leq \mathbf{Pr}\left(d^r_A(v) > \frac{d\Delta^{r-1}}{\ln^2\Delta}+\sqrt{\frac{d\Delta^{r-1}}{\ln^2\Delta}}\ln\Delta\right) < e^{-\frac{\frac{d\Delta^{r-1}}{\ln^2\Delta}\ln^2\Delta}{3\frac{d\Delta^{r-1}}{\ln^2\Delta}}} = \Delta^{-\frac{\ln\Delta}{3}} < \frac{1}{\Delta^{3r}}.\label{Av1Ineq}\end{equation}
Analogously, as the probability that $u$ belongs to $C$ equals $\frac{1}{\ln^3\Delta}$ for every $u\in N^r(v)$,
\begin{equation}
\mathbf{Pr}(A_{v,2}) \leq \mathbf{Pr}\left(d^r_C(v) > \frac{d\Delta^{r-1}}{\ln^{3}\Delta}+\sqrt{\frac{d\Delta^{r-1}}{\ln^{3}\Delta}}\ln\Delta\right) < \frac{1}{\Delta^{3r}}.\label{Av2Ineq}\end{equation}

Similarly,
\begin{equation}
\mathbf{Pr}(A_{v,3}) \leq \mathbf{Pr}\left(\left|d_A(v) - \frac{d}{\ln^2\Delta}\right| > \sqrt{\frac{d}{\ln^2\Delta}}\ln\Delta\right)
< 2e^{-\frac{\ln^{2}\Delta}{3}} < \frac{1}{\Delta^{3r}}\label{Av3Ineq}\end{equation}
and
\begin{equation}
\mathbf{Pr}(A_{v,4}) \leq \mathbf{Pr}\left(\left|d_C(v) - \frac{d}{\ln^3\Delta}\right| > \sqrt{\frac{d}{\ln^3\Delta}}\ln\Delta\right) <
\frac{1}{\Delta^{3r}}.\label{Av4Ineq}\end{equation}

Now for any $x\in [0,1]$:
$$\mathbf{Pr}(d_-(v) < X_vd-\sqrt{X_vd}\ln\Delta | X_v=x) = \mathbf{Pr}({\rm BIN}(d,x) < x d-\sqrt{x d}\ln\Delta),$$
where the probability above 
equals zero for $\sqrt{x d}\ln\Delta > x d$, while for $\sqrt{xd}\ln\Delta \leq x d$, by the Chernoff Bound,
$$\mathbf{Pr}({\rm BIN}(d,x) < x d-\sqrt{x d}\ln\Delta)
< \frac{1}{\Delta^{3r}}.$$
Hence,
\begin{equation}
\mathbf{Pr}(A_{v,5}) \leq \mathbf{Pr}(d_-(v) < X_vd-\sqrt{X_vd}\ln\Delta) \leq \int\limits_{0}^1\frac{1}{\Delta^{3r}}dx=\frac{1}{\Delta^{3r}}.\label{Av5Ineq}\end{equation}

Finally note that for $x\in[0, \frac{1}{\ln^2\Delta})$,
$$\mathbf{Pr}(d^r_-(v) > X_vd\Delta^{r-1}+\sqrt{X_vd\Delta^{r-1}}\ln\Delta \wedge v\in B | X_v=x) = 0,$$
while, analogously as above, for $x\in [\frac{1}{\ln^2\Delta},1]$:
\begin{eqnarray}
&&\mathbf{Pr}(d^r_-(v) > X_vd\Delta^{r-1}+\sqrt{X_vd\Delta^{r-1}}\ln\Delta \wedge v\in B | X_v=x) \nonumber\\
&\leq& \mathbf{Pr}({\rm BIN}(d\Delta^{r-1},x) > xd\Delta^{r-1}+\sqrt{xd\Delta^{r-1}}\ln\Delta).\nonumber
\end{eqnarray}
Since $\sqrt{xd\Delta^{r-1}}\ln\Delta \leq xd\Delta^{r-1}$ for
$x\geq \frac{1}{\ln^2\Delta}$,
by the Chernoff Bound,
$$\mathbf{Pr}({\rm BIN}(d\Delta^{r-1},x) > xd\Delta^{r-1}+\sqrt{xd\Delta^{r-1}}\ln\Delta) <
\frac{1}{\Delta^{3r}}.$$
Hence,
\begin{equation}
\mathbf{Pr}(A_{v,6}) = \mathbf{Pr}(d^r_-(v) > X_vd\Delta^{r-1}+\sqrt{X_vd\Delta^{r-1}}\ln\Delta \wedge v\in B) \leq \int\limits_{0}^1\frac{1}{\Delta^{3r}}dx=\frac{1}{\Delta^{3r}}.\label{Av6Ineq}\end{equation}

Note that each event $A_{v,i}$ is mutually independent of all other events except those $A_{u,j}$ with $u$ at distance at most $2r$ from $v$, $i,j\in\{1,\ldots,6\}$,
i.e., at most $6\Delta^{2r}+5$ events. Thus, as by~(\ref{Av1Ineq}), (\ref{Av2Ineq}), (\ref{Av3Ineq}), (\ref{Av4Ineq}), (\ref{Av5Ineq}) and~(\ref{Av6Ineq}), the probability of each such event is bounded from above by $\Delta^{-3r}$, 
by the Lov\'asz Local Lemma, with positive probability none of the events $A_{v,i}$ with $v\in V$ and $i\in\{1,\ldots,6\}$ appears. \qed
\end{pf}

Fix a vertex ordering in $V$ consistent with Lemma~\ref{MainSequencingLemma} above.

\begin{lemma}\label{SparseSpanningSubgraphInC}
There exists a set of edges $E'$ each with both ends in $C$ such that $1\leq d_{E'}(v)\leq\frac{d(v)}{\ln^{6}\Delta}$ for every $v\in C$.
\end{lemma}

\begin{pf}
For every vertex in $C=\{w_1,\ldots,w_{n'}\}$,
we choose randomly, independently and equiprobably one of its incident edges with both ends in $C$
(note that by $F_4$ from Lemma~\ref{MainSequencingLemma}, $d_C(v)\geq 1$) and denote the chosen edges by $e'_1,\ldots,e'_{n'}$, respectively (some might have been chosen twice). Set $E'=\{e'_1,\ldots,e'_{n'}\}$.
We shall argue that with  positive probability such $E'$ complies with our requirements.
For any given vertex $v\in C$ and its incident edge $e=uv$ with $u\in C$, by $F_4$,
$d_C(u)\geq \frac{1}{2}\frac{d(u)}{\ln^{3}\Delta} \geq \frac{1}{2}\ln^{5}\Delta$
(since $d(u)\geq\delta\geq \ln^8\Delta$), hence the probability that $e$ was chosen for $u$ equals at most $\frac{2}{\ln^{5}\Delta}\leq \frac{1}{6\ln^{3}\Delta}$.
As due to $F_4$, $d_C(v)\leq \frac{2d(v)}{\ln^{3}\Delta}$,
the expected number of edges incident with
$v$ chosen to $E'$ by the neighbours of $v$, denote this number by $d'_{E'}(v)$,
equals at most
$\frac{2d(v)}{\ln^{3}\Delta} \frac{1}{6\ln^{3}\Delta} = \frac{d(v)}{3\ln^{6}\Delta}$.
By the Chernoff Bound (and remarks below it), we may thus conclude that the probability of the event that $d'_{E'}(v) > \frac{d(v)}{2\ln^{6}\Delta}$, denote it by $A_v$,
is bounded from above by $e^{-\frac{d(v)}{36\ln^{6}\Delta}} \leq e^{-\frac{\ln^{2}\Delta}{36}}=\Delta^{-\frac{\ln\Delta}{36}}\leq \Delta^{-3}$.
As each such event $A_v$ is mutually independent of all other events $A_u$ except those with $u$ at distance smaller than $3$ from $v$, i.e. at most $\Delta^2$ events, by the Lov\'asz Local Lemma, with positive probability, $d'_{E'}(v) \leq \frac{d(v)}{2\ln^{6}\Delta}$ for every $v\in C$.
This however implies the existence of a desired $E'$, as $d_{E'}(v)\leq d'_{E'}(v)+ 1$ (since we must additionally only take into account the edge chosen to $E'$ by $v$ itself, which might not have been counted within $d'_{E'}(v)$) for every $v\in C$. \qed
\end{pf}

Now we shall construct a desired edge colouring applying an algorithm based on the chosen ordering of the vertices and respecting the following rules:
\begin{itemize}
  \item[(i)] We begin by attributing every edge an initial colour $Q+q$.
\end{itemize}
Then we analyze one by one subsequent vertices in the ordering and while analyzing every consecutive vertex $v$:
\begin{itemize}
  \item[(ii)] we allow adding or subtracting $Q$ to the colour of every backward edge of $v$;
  \item[(iii)] we allow adding any integer in $\{0,\ldots,q\}$ to the colour of every forward edge of $v$.
\end{itemize}
Note that after introducing such changes, for the obtained final colour $c(e)$ of every edge $e\in E$ we shall have:
\begin{equation}\label{weight_bounds}
q\leq c(e)\leq 2Q+2q,
\end{equation}
as desired.
\begin{itemize}
  \item[(iv)] Special rules shall be used in the final part of the construction, but these shall still be consistent with the inequalities from (\ref{weight_bounds}) above.
\end{itemize}
\begin{remark}\label{DegreeRemark}
Note that by the bounds from (\ref{weight_bounds}) above, since $\frac{2Q+2q}{q}<5\ln\Delta$, any $r$-neighbours $u,v$ with
$$d(u) \geq d(v) 5\ln\Delta$$
shall certainly be sum-distinguished in $G$ at the end of our construction.
\end{remark}

Let us define the following family of $2$-element sets of integers:
$$\mathcal{P}=\{\{p,p+Q\}|p\in\{0,\ldots,Q-1\} {\rm ~mod~} {2Q}\}.$$
Note that the sets in this family are pairwise disjoint.
Every consecutive vertex $v$ in the sequence shall have assigned one of such sets $W_v\in \mathcal{P}$ (the moment it is analyzed),
and ever since such assignment, the sum at $v$ shall always be required to belong to $W_v$.

Let us begin analyzing the consecutive elements of the sequence starting from the first one, thus firs we consider subsequent vertices in $A$.

Note that every vertex $v\in A$ of degree $d$, by $F_1$ (from Lemma~\ref{MainSequencingLemma}), has at most
$\frac{2d\Delta^{r-1}}{\ln^2\Delta}$ backward $r$-neighbours, from which it has to be distinguished.
At the same time,
by $F_4$ it has at least $\frac{d}{2\ln^{3}\Delta} > \frac{Q}{q}$ forward neighbours (edges).
Therefore, using admissible alterations from (ii) and (iii) on the backward edges (to the colours of which we may either add or subtract $Q$ so that the sums of their initial ends remain in their corresponding already assigned 2-element sets 
from $\mathcal{P}$) and forward edges incident with $v$,
we may obtain at least
$d q+1$
consecutive integer sums at $v$, among which there are elements (not necessarily both) from at least
$\frac{d q}{2}$
pairs from $\mathcal{P}$.
Moreover, the elements in at least $\lfloor\frac{1}{3}\frac{d q}{2}\rfloor \geq \frac{d\Delta^{r-1}}{7\ln\Delta} > 2\frac{d\Delta^{r-1}}{\ln^2\Delta}$
of such pairs are congruent to $0$ modulo $3$ (note that by the choice of 
$Q$, the two elements in every pair in $\mathcal{P}$ are congruent modulo $3$).
We thus may perform admissible alterations on the edges incident with $v$ so that afterwards it has a sum belonging to some pair in $\mathcal{P}$ with elements congruent to $0$ modulo $3$ which is disjoint with 
all $W_u$ associated with backward $r$-neighbours $u$ of $v$ (this way we shall among others guarantee the distinction between sums of all vertices in $A$).
 We also set this pair as $W_v$ and continue in the same manner with all vertices in $A$.

Suppose now that $v\in B$ has degree $d$, and thus far all our rules and requirements have been fulfilled.
Analogously as above, by $F_4$, $F_5$ and admissible operations (ii), (iii), we have
at least
$Q(X_v d-\sqrt{X_v d}\ln\Delta)+q[d-(X_v d-\sqrt{X_v d}\ln\Delta)]\geq 2\Delta^{r-1}(X_v d-\sqrt{X_v d}\ln\Delta) + d \frac{\Delta^{r-1}}{\ln\Delta}$
consecutive integer sums available for $v$, including
elements from at least $\Delta^{r-1}(X_v d-\sqrt{X_v d}\ln\Delta) + d \frac{\Delta^{r-1}}{2\ln\Delta}$
pairs from $\mathcal{P}$,
while by $F_6$, there are at most
$X_v d \Delta^{r-1}+\sqrt{X_vd\Delta^{r-1}}\ln\Delta$
backward $r$-neighbours of $v$, where
$$\Delta^{r-1}(X_v d-\sqrt{X_v d}\ln\Delta) + d \frac{\Delta^{r-1}}{2\ln\Delta} > X_v d \Delta^{r-1}+\sqrt{X_vd\Delta^{r-1}}\ln\Delta,$$
hence we may set the sum of $v$ via admissible alterations from (ii) and (iii) so that it belongs to a set in $\mathcal{P}$ disjoint with $2$-element sets associated with all its backward $r$-neighbours, and fix this set as $W_v$.
We apply such greedy algorithm to all consecutive vertices in $B$.

Now,
before we shall continue with the vertices in $C$,
we first subsequently analyze all edges $e_1,e_2,\ldots,e_{n''}$ in $E'$ guaranteed by Lemma~\ref{SparseSpanningSubgraphInC}
(which induce a spanning subgraph in $G[C]$),
and choose additions to their colours in the range $0,1,2$ greedily so that afterwards the sum at every vertex $v\in C$ is not congruent to $0$ modulo $3$.
Then we randomly, independently and equiprobably subtract from the colour of every edge in $E'$ an integer in $\{0,\ldots,Q-1\}$
divisible by $3$; these shall be the final alterations of colours of the edges in $E'$. (Note that afterwards, $q\leq c(e)\leq 2Q+2q$ for every such edge $e$, hence (\ref{weight_bounds}) is fulfilled for these.)
Denote the (temporary) sum obtained for every $u\in C$ by $w'(u)$ (and note that $w'(u)\not\equiv 0 {\rm ~mod~} 3$).

\begin{lemma}\label{LemmaDistributionInC}
With positive probability,
for each vertex $v\in C$ and every integer $t\in [0,Q-1]$ which is not congruent to $0$ modulo $3$,
the number of vertices $u$ in $N_C^r(v)$ with $(5\ln\Delta)^{-1} d(v)\leq d(u) \leq d(v) 5\ln\Delta$ and $w'(u)\equiv t {\rm ~mod~} Q$ is upper-bounded by $5\frac{d(v)}{\ln^{3}\Delta}$.
\end{lemma}

\begin{pf}
For each vertex $v\in C$ of degree $d$ in $G$ and every integer $t\in [0,Q-1]$ which is not congruent to $0$ modulo $3$,
let $X_{v,t}$ denote the number of vertices $u$ in $N_C^r(v)$ with $w'(u)\equiv t {\rm ~mod~} Q$
and $(5\ln\Delta)^{-1} d\leq d(u) \leq d 5\ln\Delta$. As for every
$u\in N_C^r(v)$ with $(5\ln\Delta)^{-1} d\leq d(u) \leq d 5\ln\Delta$,
$\mathbf{Pr}(w'(u)\equiv t {\rm ~mod~} Q) \leq \frac{3}{Q}$
(what can be easily proved by means of the total probability via analysis of the possible $\frac{Q}{3}$ choices of subtractions for the last edge in $\{e_1,\ldots,e_{n''}\}=E'$ incident with $u$, at most one of which assures
$w'(u)\equiv t {\rm ~mod~} Q$ regardless of any fixed choices for the remaining edges),
by $F_2$ we thus obtain that
$\mathbf{E}(X_{v,t})\leq \frac{3}{Q} \frac{2d\Delta^{r-1}}{\ln^{3}\Delta}
\leq 3 \frac{d}{\ln^{3}\Delta}$.

Note that a change of choice for any edge in $E'$ may influence $X_{v,t}$ by at most $2$.
Moreover, for any $s$, the fact that $X_{v,t}\geq s$ can be certified by the outcomes of at most
$s\cdot \frac{5d}{\ln^{5}\Delta}$ 
trials, i.e., choices committed on the edges in $E'$ incident with some $s$ $r$-neighbours $u$ of $v$ in $C$
with $(5\ln\Delta)^{-1} d\leq d(u) \leq d 5\ln\Delta$, each of which has at most $\frac{d 5\ln\Delta}{\ln^{6}\Delta}=\frac{5d}{\ln^{5}\Delta}$ incident edges in $E'$ by Lemma~\ref{SparseSpanningSubgraphInC}.
Thus
by Talagrand's Inequality (and comments below it),
\begin{eqnarray}
&&\mathbf{Pr}\left(X_{v,t} > 5 \frac{d}{\ln^{3}\Delta}\right)\nonumber\\ 
&\leq&
\mathbf{Pr}\left(X_{v,t} > 3 \frac{d}{\ln^{3}\Delta} + \frac{d}{\ln^{3}\Delta}+ 20\cdot 2\sqrt{\frac{5d}{\ln^{5}\Delta} 3 \frac{d}{\ln^{3}\Delta}}
+64\cdot 2^2 \frac{5d}{\ln^{5}\Delta}\right)\nonumber\\
&<& 4e^{-\frac{\left(\frac{d}{\ln^{3}\Delta}\right)^2}{8\cdot 2^2 \frac{5d}{\ln^{5}\Delta} 
\left(3 \frac{d}{\ln^{3}\Delta} + \frac{d}{\ln^{3}\Delta}\right)}}<\frac{1}{\Delta^{5r}}.\label{PrXvtBound}
\end{eqnarray}
As any event that $X_{v,t} > 5 \frac{d}{\ln^{3}\Delta}$ is mutually independent of all other events of the form $X_{v',t'} > 5 \frac{d(v')}{\ln^{3}\Delta}$ with $d(v,v')>2r+1$, i.e., all except at most $\Delta^{2r+1}\cdot \frac{2Q}{3}<\Delta^{4r}$ such events,
by the Lov\'asz Local Lemma and (\ref{PrXvtBound}) we thus obtain the thesis.
\qed
\end{pf}

We fix any subtractions from the colours of the edges in $E'$ consistent with the thesis of Lemma~\ref{LemmaDistributionInC}.
Then,
as by $F_4$ every vertex in $B$ has a neighbour in $C$, we subtract $Q$ if necessary (or do nothing) from the colour of one such edge for every vertex in $B$ so that the weight for every vertex $v\in B$ is set on the smaller element of its associated two-element list $W_v$.
(Note that prior to these changes, every such edge had its colour between $Q+q$ and $Q+2q$, as it has not been analyzed as a backward edge yet, hence
(\ref{weight_bounds}) shall hold for this edge after any of the described changes).

Note that by our construction, the sums at the vertices in $C$ are not congruent to $0$ modulo $3$, contrary to the sums at vertices in $A$, hence vertices in $A$ and $C$ are distinguished from each other.
This shall not change till the end of the construction, as while analyzing the consecutive vertices in $C$ we shall only allow adding or subtracting $Q$ on the edges between $A$ and $C$ so that the sums of the vertices in $A$ remained in their associated $2$-element sets.
While performing these changes we need only guarantee the distinction between vertices in $C$ (and creating no conflicts between vertices in $C$ and $B$). Note that via these admitted operations, for every $v\in C$ with $d(v)=d$ and $w'(v)\equiv t {\rm ~mod~} Q$,
by $F_3$ we may obtain at least $\frac{d}{2\ln^{2}\Delta}+1$ consecutive sums congruent to $t {\rm ~mod~} Q$ at $v$,
among which we have at least $\frac{d}{4\ln^{2}\Delta}$ available options that are not used as a sum of any vertex in $B$ (nor obviously in $A$), since these are all fixed on the lower positions from their associated lists.
As by Lemma~\ref{LemmaDistributionInC} (and Remark~\ref{DegreeRemark}) above we need to only distinguish $v$ from its at most
$5\frac{d}{\ln^{3}\Delta} < \frac{d}{4\ln^{2}\Delta}$
$r$-neighbours $u$ in $C$ (with $w'(u)\equiv t {\rm ~mod~} Q$), we have at least one available choice for the sum of $v$ consistent with this goal
(note that this time the sum at $v\in C$ may belong to a list $W_{v'}$ associated with some vertex $v'\in B$).
After analyzing all consecutive vertices in $C$ we thus obtain an edge colouring of $G$ with colours in $[q,2Q+2q]$ and without any conflicts between $r$-neighbours in $G$.
\qed

\section{Remarks}
We note that the lower bound of $\ln^8\Delta$ for $\delta$ in Theorem~\ref{przybylo_main_th}
was chosen, but certainly not optimized,
for the sake of clarity of presentation.
Nevertheless, our approach does not seem to allow to remove any poly-logarithmic (in $\Delta$) lower bound on the minimum degree of a graph.

The upper bound from Theorem~\ref{przybylo_main_th} can on the other hand
be significantly improved in terms of the magnitude of the second order term, i.e. $\frac{\Delta^{r-1}}{\ln\Delta}$,
in the case of graphs with relatively large minimum degree, e.g. for regular graphs, using even a slightly simplified version of the algorithm presented in Section~\ref{SectionWithPrzybyloProof}.

We conclude by posing a conjecture which to our believes 
expresses
a true asymptotically optimal upper bound for the investigated parameters.
\begin{conjecture}\label{przybylo_main_con}
For every integer $r\geq 2$ and each
graph $G$ with maximum degree $\Delta$ and without an isolated edge,
$$s_r(G) \leq (1+o(1))\Delta^{r-1}.$$
\end{conjecture}
We also refer a reader to~\cite{Przybylo_distant_total_probabil} to see an improvement of a similar probabilistic flavor for the upper bound from~\cite{Przybylo_distant} on the correspondent of $s_r(G)$ concerning the case of total colourings.

\end{document}